%
\parindent=0pt\parskip=0pt\magnification=\magstep0


{\bf Errata, updates of the References, etc.,
(as of Jan. 21, 2004) for the book:}

\vskip .5cm

Basic Hypergeometric Series, by George Gasper and Mizan Rahman,
Encyclopedia of Mathematics and its Applications, Vol. 35, Cambridge University
Press, Cambridge-New York, 1990. xx+287$\,$pp. \hfil\break
ISBN 0-521-35049-2

\vskip .3cm

{\bf Note: An expanded updated second edition of this book is due to 
be published during 2004.
It will contain three additional chapters entitled:}
\vskip .1cm
{\bf\ 9.  Linear and Bilinear Generating Functions for Basic 
Orthogonal Polynomials}

{\bf 10. $q$-Series in Two or More Variables}

{\bf 11.  Elliptic, Modular, and Theta Hypergeometric Series}

\vskip .5cm
Errata added after October 16, 1995 are placed at the end of
the following list.


\vskip .2cm

p. xv, line 20: Replace ``seems''  by ``seem''.

p. 5: Replace `` $(1 - q^n)$ '' in (1.2.26) by  `` $(1 - q^{n+1})$ ''.

p. 5, third paragraph: Replace each `` $b_s|$ '' by `` $b_s q|$ '' .

p. 5: Delete  `` $/a$ ''  from the right side of eq. (1.2.27).

p. 11, line 4 up: Replace  ``and''  by  ``by''.

p. 12:   Insert the factor `` $z^n$ ''  in the numerator on the
left side of eq. (1.6.3).

p. 12, line 7 up: Replace `` $n = 0$ ''  by `` $k = 0$ ''.

p. 15, line 15: Replace `` $q^{k \choose 2}$ '' by `` $q^{-{ k \choose 2}}$ ''.

p. 15, line 16 (the last line of (1.9.4): Replace `` $; $ '' by `` $; q, $ ''.

p. 17, third line below (1.10.5):  Replace ``analogue'' by  ``$q$-analogue''.

p. 19, line above eq. (1.11.8): Replace ``integal'' by ``integral'' .

p. 24, Ex. 1.19 (ii): Replace `` $(a^2 q, b^2 q; q^2)$ '' by
  `` $(a^2 q, b^2 q; q^2)_\infty$ ''.

p. 24, Ex. 1.21: Replace `` $1 - q^{n+x}$ '' by `` $(1 - q^{n+x})^2$ ''.

p. 25, Ex. 1.24, line 7:  Replace `` $J_\nu ^{(1)}(x;q)_\infty$ ''  by
`` $J_\nu ^{(1)}(x;q)$ ''.

p. 26, line 2:   Replace  `` $x/2t$ '' by `` $-x/2t$ ''.

p. 26, Ex. 1.29, line 1: Replace ``ultraspherial'' by ``ultraspherical''.

p. 26, Ex. 1.29, line 6: Replace `` $e^{-in\theta}$ ''
by `` $\beta^{-n} e^{-in\theta}$ ''.

p. 27, Ex. 1.32: Place a period at the end of the first sentence.

p. 28, Ex. 1.35: Replace ``associate law''  by  ``associative law'' .

p. 28, line 3 up:  Replace  ``derived''  by ``derive''.

p. 36, line 4: Replace ``(2.1.5)'' by ``(2.1.6)''.

p. 36, line 6: Replace ``(2.1.5)'' by ``(2.1.6)''.

p. 36, line 10: Replace  ``q-factor'' by ``$q$-factor'' (i.e, with a math.
italic $q$).

p. 36, line below eq. (2.7.1): Replace ``(2.1.6)''  by  ``(2.1.7)''.

p. 36, line above eq. (2.7.2): Replace ``(2.1.4)'' by ``(2.1.5)''.

p. 37, line 1: Replace `` $a q^{n+1}$ '' by  `` $-a q^{n+1}$ ''.

p. 37, Eq. (2.7.7): Replace `` $ ; q^5)$ ''  by  `` $; q^5)_\infty$ '' .

p. 41, line 4 up: Replace ``respectivey''  by  ``respectively''.

p. 43: A sub infinity is missing from the right side of the numerator
of the first fraction on the right side of (2.10.12).

p. 47, line 4: Replace `` $a/bq^n$ and $a$ ''  by `` $\lambda e/a$ 
and $a/bq^n$ ''.

p. 50, Ex. 2.2:  Replace  `` $(aq, bqt; $ ''  on the right side by
`` $(aq, bt;  $ ''      (see eq. (2.3) in Gasper and Rahman [1983a]).

p. 51, Ex. 2.13 (i): In the $_{10}W_9$  replace `` ,; ''  by `` ; '' 
and replace
the argument `` $\lambda/a$ ''  by  `` $q \lambda /a $ ''.

p. 52, Ex. 2.14 (ii): After the equation add ``, where $\lambda = a^2 
q/bcd$ '' .

p. 54, Ex. 2.21: Replace `` $({aq \over \lambda})^n$ '' by
`` $({a^2q^2 \over efgh})^n$ ''.

p. 56, Ex. 2.30: Replace `` where $a^3 q^2 = bdefgh$ '' by
`` where $a^3 q^2 = bcdefgh$ '' .

p. 57, Notes 2, line 5: Add  ``Adiga {\it et al.} [1985],''.

p. 59:  Replace the  `` $q$ ''  in eq. (3.1.1)  by   `` 1 ''.

p. 62, line below eq. (3.2.9):  Replace ``min'' by ``max'' .

p. 66, line 4 up:  Replace `` $x/aq$ ''  by `` $x/bq$ ''.

p. 68:  Replace the third  `` ; ''  on the left side of eq. (3.5.5) 
by a `` , ''.

p. 73, eq. (3.6.19): Replace  `` $q)_{n-1}$ '' by `` $p)_{n-1}$ '' .

p. 74:  Replace the second  `` $n = 0$ '' on the right side of eq. (3.7.2)
by  `` $k = 0$ '' .

p. 78, eq. (3.8.1):  Replace the factor `` $p^k$ '' on the right side of the
summand by `` $p^{2 k}$ '' .

p. 85, line below eq. (3.9.2): Replace   ``in (3.9.2)'' by  ``in (3.9.1)''.

p. 88, eq. (3.10.5):   Replace the `` $q^{n+1}$ '' in the argument of 
the series
by `` $q^{n-1}$ '' .  The right side of the equation may be simplified to
`` $(-aq, aq^2/w, w/aq; \, q)_n / (-q, aq/w, w; \, q)_n$ '' .

p. 88, eq. (3.10.7):   Replace `` $(-a q^2,$ '' on the right side of 
the equation
by `` $(-a q,$ '' .

p. 88, eq. (3.10.9): The right side of the equation may be simplified to
`` $(w/aq, -a^{1/2}, aq^2/w; \, q)_n \over (w, -a^{-1/2},aq/w; \, q)_n$ '' .

p. 89, eq. (3.10.14): The fraction on the right side of the equation may be
simplified to  `` $(-aq, aq^2/w, w/aq; \, q)_n \over  (-q, aq/w, w;\, q)_n$''.

p. 89:  In the numerator of the display at the bottom of the page, replace
``$ -\lambda q^{n+1}$ '' by `` $-\lambda q^{n+1}/a$ '' .

p. 90, Ex. 3.2 (ii): Replace `` $(a^2 z^2;q)_\infty$ '' by
`` $(a^2 z^2;q^2)_\infty$ '' .

p. 91, Ex. 3.4: Replace the last `` ${}_4 \phi_3$ '' by `` ${}_4 \phi_2$ ''
and delete the `` $, 0$ '' from its denominator parameters.

p. 93, Ex. 3.10, line 3: Replace `` $(\lambda q,$ '' by
`` $(\lambda/a, \lambda q,$ '' and replace `` $(\lambda q/a,$ '' by
`` $(\lambda, \lambda q/a,$ '' .

p. 95,  Ex. 3.21: In the numerator on the right side of the equation replace
`` $(A, AQ/B; Q)_n$ '' by \break `` $(Q, AQ/B; Q)_n$ '' .

p. 95, Ex. 3.21, line 5: Replace the upper limit `` $\infty$ '' on the
summation symbol by  `` $n$ ''.

p. 96, Ex. 3.21:  In the left side of the last equation, replace
`` $q/BCq^n$ '' by `` $1/BCq^n$ '' .

p. 97, Ex. 3.25: Replace the `` $a^3$ '' in the denominator of the 
first line by
`` $a^2$ ''.

p. 98, Ex. 3.29 (ii): Delete the second `` $(ac/q^2;q)_k$ ''.

p. 99:  Delete the period from the end of the last line.

p. 100, line 2 up:  Replace  ``Ex. 3.16'' by  ``Ex. 3.11''.

p. 101, line 1 up: Delete ``and Rother'' .

p. 104, line 8:  Replace ``anlaytic''  by ``analytic''  .

p. 118:  Insert a left parenthesis onto the left side of the
numerator on the right side of eq. (4.11.2).

p. 122, Ex. 4.5:  Replace `` ${}_4 \phi _2$ '' by `` ${}_4 \phi _3$ '' .

p. 123, line 9: Replace `` ${\rm a^3 c^5}$ '' by `` $a^3 c^5$ ''.

p. 124, first line of Notes 4: Replace ``[1989b]'' by ``[1989a,b]'' .

p. 134, Ex. 5.3: Replace `` In (5.4.1) set $b = a, d=c, e=-1$ ''
by `` In (5.3.1) set $b = a, c = d = e = -1$ ''.

p. 135, Ex. 5.10: In the denominator on the right hand side of the 
equation replace
  `` $(c,$ '' by `` $(e,$ '' .  Insert a period at the end of the equation.

p. 136:  Insert a period at the end of the equation in Ex. 5.15 .

p. 137, Ex. 5.18 (iii):  Add `` when $bcde = q^{n+1},$ '' .

p. 137, Ex. 5.18 (iv):  Add `` when $bcde = q^{n+3},$ '' .

p. 147, eq. (6.4.11): In the numerator of the first ${}_{10} \phi_9$ 
on the right
  side of the equation replace `` $\nu_1, \sqrt \nu_1$ '' by `` 
$\nu_1,  q \sqrt \nu_1$ '' .

p. 158, line 2 up: Replace ``6.14'' by  ``6.13''.

p. 167, eq. (7.3.13): Replace the factor $(-ac)^n$ by $(-acq^2)^{-n}$ .

p. 170, line 6 up: Replace ``[1986a,b, 1987]'' by  ``[1986a]''.

p. 176, line 2 up:  Replace ``q-analogues'' by  ``$q$-analogues''.

p. 177:  On the right side of eq. (7.5.31) replace
`` $ -n(2 \alpha +1)/4 $ ''  by  `` $ -n(2 \alpha +1)/2 $ '' .

p. 177: On the right side of eq. (7.5.32) replace
`` $ q^{\beta +1}  -q^{\alpha + \beta +1} $ '' by
`` $ q^{\beta +1},  -q^{\alpha + \beta +1} $ '' .

p. 180, Eq. (7.6.5): Replace each of the two `` $(ad;q)_j$ '' by
`` $(ad,bd,cd;q)_j$ '' and replace `` $(abcdq;q)_j$ '' by `` $(abcd;q)_j$ ''.

p. 181, formula (7.6.9): Replace `` $(bc,bd,cd;q)_n$ ''
by `` $(q,bc,bd,cd;q)_n$ ''.

p. 183, line 7 up: Insert `` $\sqrt {1 - x^2}$ ''  as a factor in
front of `` $w(x;a,b,c,d|q)$ ''.

p. 183, eq. (7.7.16): Delete  `` $\sqrt {1 - x^2}$ '' (three of them).

p. 184, eq. (7.7.17): Delete  `` $\sqrt {1 - x^2}$ '' (two of them).

p. 186, Ex. 7.11: Insert the missing right parenthesis onto the last
exponent in the fourth line.

p. 187, Ex. 7.13, line 5: Replace the multiple of the delta function on the
right side of the orthogonality relation by its reciprocal.

p. 191, Ex. 7.32 (ii):  Insert the factor `` $q^{(\alpha - \beta) n/2}$ ''
after the  `` $(-1)^n$ '' factor on the right side of the equation.

p. 192:  Delete the four minus signs from the left side of the
equation at the top of page 192.

p. 194, Ex. 7.43, line 2: In the $_1\phi_1$ series replace `` $-x$ ''
by `` $-x(1-q)$ ''.

p. 197, line 2:  Replace ``lime'' by ``line''.

p. 202, line 2: Replace the last `` $-m$ '' in the denominator
parameters in (8.2.12) by `` $m$ ''.

p. 207, line 2 below (8.5.1): Replace ``appers''  by  ``appears''.

p. 212, eq. (8.7.13):  Insert the factor `` $(q; q)_s$ '' into the
denominator on the right side of the equation.

p. 223, line 12:  Replace  ``partion''  by  ``partition''.

p. 226, line 4:  Replace ``Ex. (5.1) and (5.4)'' by ``Ex. 5.1 and 5.4'' .

p. 232, second $_5\phi_4$  in Ex. 8.17: Replace `` $z^2 b^2/z^2$ ''
by `` $a^2 b^2/z^2$ ''.

p. 232, Ex. 8.18: Add the restriction that $0 < q < 1$.

p. 232, line 11 up: Insert   ``de Branges and Trutt [1978],''.

p. 234, (I.31): Replace `` $,q$ '' by  `` $; q$ ''.

p. 237:  Insert the factor `` $q^{n(n+1)/2}$ '' in front of the right hand
side of eq. (II.19).

p. 239, line 6 up: Replace  ``bilaterial'' by   ``bilateral''.

p. 240, eq. (II.37): Replace the first factor `` $(p; p)_n$ '' in
the denominator by `` $(p; p)_k$ '' .

p. 241, Eq. (III.5): Add the restriction that `` $ a = q^{-n}$, where  $n$  is
a nonnegative integer''. (See \break Ex. 1.15 (i).)

p. 241, (III.8):  Replace `` $q/c$ '' by  `` $q/z$ ''.

p. 243, (III.23): Replace the second `` $a^{1\over 2}$ '' in the denominator of
the first $_8\phi_7$ by `` $-a^{1\over 2}$ ''.

p. 253: In Askey [1975] delete ``, ed.'' and replace ``Applies'' by ``Applied".
Askey [1989a] is published in Ramanujan International Symposium on Analysis
(N.K. Thakare, ed.),
Macmillam India, Delhi, pp. 1--83.  Russian translation, Usp.\ Math.\
Nauk.\ 45(1) (1990), pp. 33--76.
Also Russian Math.\ Surveys 45:1 (1990), pp. 37--86.

Askey [1989b] is published in $q$-Series and Partitions (D. Stanton, ed.),
Springer, Berlin and New York, pp. 151--158.

Askey [1989c] is published in Part 2, pp. 299--321.

Askey [1989d] is published in Number Theory and Related Topics, Ramanujan
Colloquium, Bombay (1988), Oxford Univ.\ Press, pp. 1--12.

Askey [1989e] is published in Number Theory (K. Alladi, ed.), Lecture
Notes in Math.\ 1395, Springer, Berlin and New York, pp. 84--121.

p. 254: In Askey [1989f] replace ``to appear'' by  ``pp. 121--129''.

Askey [1989g] is published in Computers in Mathematics (D. Chudnovsky
and R. Jenks, eds.),

M. Dekker, New York and Basel (1990), pp. 35--82.

In Askey and Wilson [1979], replace ``generalized'' by ``generalize''.

p. 256: Change ``de Branges, L. (1978)'' to ``de Branges, L. and Trutt, D.
(1978)'' and add (after New York) ``, pp. 1--24''.

p. 261: In Gasper [1985] delete the hyphen that is between ultra and spherical.

Gasper [1989d] is published in $q$-Series and Partitions (D. Stanton, ed.),
IMA Vol.\ Math.\ Appl.\ 18 (1989), Springer, Berlin and New York, pp. 15--34.

Gasper [1989e] is published in Orthogonal Polynomials: Theory and Practice
(P. Nevai, ed.), Kluwer Academic Publishers, Boston, 1990, pp. 163--179.

p. 262: Gasper and Rahman [1989a] is published in Canad.\ J.\ Math.\
42 (1990), 1--27.

In  Gasper and Rahman [1989b] change ``to appear'' to ``20, 1270--1282''.

p. 268: Koornwinder [1989a] is published in J.\ Math.\ Anal.\ Appl.\
148(1990), 44--54.

Koornwinder [1989b] is published in SIAM J.\ Math.\ Anal.\ 22(1991), 295--311.

In Macdonald [1972] replace ``$h$-function''  by  ``$\eta$-function''.

p. 269: Milne [1989c] is published in {\it Theta Functions,  Bowdoin 1987},
Part 2 (L. Ehrenpreis
and R. Gunning, eds.), Amer.\ Math.\ Soc., Providence, R.I., pp. 323--359.

p. 270: Replace ``Paule, P. and Rother, W. (1985). ... 1-37.'' by
``Paule, P. (1985). Ein neuer Weg zur {\it q}-Lagrange Inversion,
Bayreuther Math.\ Schriften 18, 1-37.''

p. 271:  In Rahman [1989a] change ``to appear'' to ``Proc.\ Amer.\ Math.\
Soc.\ 107, 373--381''.

Rahman [1989b] is published in {\it International Symposium on Analysis}
(N.K. Thakare, ed.), Macmillan India, pp. 117--137.

In Rahman [1989c] change the title and location to ``A note on
the orthogonality of Jackson's $q$-Bessel
functions, Canad.\ Math.\ Bull.\ 32(3), 369--376''.

In Rahman [1989e] change ``to appear'' to ``Utilitas Mathematica 36, 
161--172''.

p. 274: In Stanton [1989] delete ``{\it Workshop on}'' and
replace ``to appear'' by ``pp. 139--149''.

p. 275: In Wilson [1978] change ``polynomials'' to ``functions'' .

p. 279: Delete ``Rother, W., 101, 270'' .

p. 280, below Trjitzinsky: add  ``Trutt, D., 232, 256''.

p. 282: Add the line `` $w(x;a,b,c,d|q) \quad 174$ ''  below  the line
`` $w(x;a,b,c,d) \quad  143$ ''.

p. 282, line 3 up: Replace `` $U_n^{(\alpha)}$ '' by `` $U_n^{(a)}$ '' .

\vskip .3cm

{\bf The following errata were added (in the order listed) after 
October 16, 1995:}
\vskip .2cm

p. 88, eq. (3.10.8):   Replace `` $, a/q^2;$ '' on the right side of the
equation by `` $, a/q b^2;$ '' .

pp. 57, 232 and 256:  Delete the reference to Berndt and Joshi [1983]

p. 128,  eq. (5.3.1):   Replace the `` ; '' on the right side of the
equation by a `` , '' .

p. 129, line above eq. (5.3.4):  Replace `` very-well-poised '' by
`` well-poised '' .

p. 260: Garvan [1989] is published in SIAM J.\ Math.\ Anal.\ 21 
(1990), 803--821.

Garvan and Gonnet [1989] is published in J. Symbolic Comput.\ 14 
(1992), 141--177.

Garvan and Stanton [1989] is published in Math.\ Comp.\ 55 (1990), 299--311.

p. 263: Gustafson [1989a] is published in SIAM J. Math.\ Anal.\ 21 
(1990), 510--522.

Gustafson [1989b] is published in Ramanujan International Symposium on Analysis
(Pune, 1987), 185--224, Macmillan of India, New Delhi, 1989.

  Habsieger [1989] is published with the title ``Laurent Une $q$-int\'egrale
de Selberg et Askey'' (French) [A $q$-integral of Selberg and Askey] in
SIAM J. Math.\ Anal.\ 19 (1988), 1475--1489.

p. 264: Ismail, Merkes, and Styer [1989] is published in Complex Variables
14 [1990], 77--84.

p. 266: Joichi and Stanton [1989] is published in Discrete Math.\ 73 (1989),
261--271.

p. 267: Kadell [1989c] is published in
  Mem.\ Amer.\ Math.\ Soc.\ 108 (1994), no. 516.

Kalnins and Miller [1989a] is published in Constructive Function Theory---86
Conference (Edmonton, AB, 1986), Rocky Mountain J. Math.\ 19 (1989), 223--230.

Kalnins and Miller [1989b] is published in  SIAM J. Math.\ Anal.\ 19
(1988),1216--1231.

In Kirillov and Reshetikhin [1988] change the title and location to 
``Representations
of the algebra ${U}\sb q({\rm sl}(2)),\;$ $q$-orthogonal polynomials 
and invariants of
links. Infinite-dimensional Lie algebras and groups 
(Luminy-Marseille, 1988), 285--339,
  Adv.\ Ser.\ Math.\ Phys., 7, World Sci.\ Publishing, Teaneck, NJ, 1989.''

p. 268: Krattenthaler [1989a] is published in Canad.\ J. Math.\ 41 
(1989), 743--768.

Krattenthaler [1989b] is published in Proc.\ Amer.\ Math.\ Soc.\
124 (1996), 47--59.

Krattenthaler [1989c] is an unpublished manuscript; a summary of the results
in it appeared in: Einige quadratische, kubische und quartische 
Summenformeln f\"ur
$q$-hypergeometrische Reihen, Anzeiger \"OAW math.-naturw. Klasse 126 
(1989), 9--10.

The second edition of Macdonald [1979] with contributions by A. 
Zelevinsky is published
  in Oxford Mathematical Monographs, Oxford University Press, Oxford 
and New York, 1995.

p. 269: Macdonald [1989] is published in Orthogonal polynomials 
(Columbus, OH, 1989),
311--318, NATO Adv.\ Sci.\ Inst.\ Ser.\ C Math.\ Phys.\ Sci., 294, 
Kluwer Acad.\
Publ.,  Dordrecht, 1990.

Masuda, Mimachi, Nakagami, Noumi, and Ueno [1989] is published in
J. Funct.\ Anal.\ 99 (1991), 357--386.

Milne [1989a] is published in Adv.\ in Math.
  99 (1993), 162--237.

Milne [1989b] is published in Adv.\ in Math.\ 108 (1994), 1--76.

Mimachi [1989] is published in Nagoya Math.\ J. 116 (1989), 149--161.

p. 270: O'Hara [1989] is published in J. Combin.\ Theory Ser.\ A 53 
(1990),  29--52.

p. 271: Rahman[1989d] is published in Canad. J.\ Math.\ 45 (1993), 394--411.

In Yang [1989] change the title and location to ``Partitions, Inversions,
and $q$-Algebras, SIAM Proceedings on Graph Theory, Combinatorics, 
Algorithms, and
Applications, ed. by F.R.K. Chung, R.L. Graham, and D.F. Hsu, (1991), 
pp. 625--635.''

p. 276: Zeilberger [1989a] is published in Discrete Math.\ 79 (1990), 313--322.

Zeilberger [1989b] is published in Amer.\ Math.\ Monthly 96 (1989), 590--602.

Zeilberger [1989c] is published in J. Comput.\ Appl.\ Math.\ 32 
(1990), 321--368.

Zeilberger [1989d] is published in  $q$-Series and Partitions (D. 
Stanton, ed.),
IMA Vol. Math. Appl. 18 (1989), Springer, Berlin and New York, pp. 35--44.

p. 278: Delete ``Joshi, P.T., 57, 232, 256'' .

p. 122, Ex 4.5: Replace `` $(q^{1 \over 2},q^{1 \over 2};q)_{\infty}$ '' by
`` ${(ab, ac, bc,q^{1 \over 2},q^{1 \over 2};q)_{\infty}}$ '' .

p. 52, Ex. 2.14 (ii): In the third numerator parameter of the $_{12}\phi_{11}$
replace `` $\lambda^{-{{1}\over{2}}}$ '' by `` $\lambda^{{1}\over{2}}$ ''.

p. 90, eq. (3.10.16): In the last denominator parameter of the $_{10}\phi_9$
replace `` $q^{2n+1}$ '' by `` $q^{2n+2}$ '' .

p. 93, Ex. 3.13, first identity: Replace the argument of the $_{10}W_9$
by `` $a^2 q^3/ \lambda$ '' .

p. 93, Ex. 3.13, second identity: In the denominator of the coefficient on
the right side of the equation replace `` $(-q, \lambda/a^2; q)_n$ ''
by `` $(-q, \lambda/a^2 q; q)_n$ '' .

p. 94, Ex. 3.14:  On the right side of the equation replace
`` $({2 \over x})^\nu$ '' by `` $({x \over 2})^\nu$ '' .

p. 235:   In the {\bf q-Binomial coefficient} section it is assumed that
$n$ and $k$ are nonnegative integers.

pp. 247--248, eq. (III.40):  In the argument of the ${}_{10}\psi_{10}$
and the ${}_{10}\phi_9$ replace `` $q^2$ ''  by `` $q^3$ ''.  In the
denominator of the coefficient on the right side of the equation replace
`` $gk, h/g$ '' by `` $gk, k/g,  h/g$ '' .

p. 147, eq. (6.4.10): In the coefficient on the right side of the equation
replace `` $\tau d, \sigma/\tau $ ''  by `` $\tau d, \tau f, \sigma/\tau $''
and replace `` $ \sigma/d, \tau^2  $ ''  by
`` $ \sigma/d, \sigma/f, \tau^2  $ '' .

p. 27, Ex. 1.33, line 4: Insert the factor
  `` $(-1)^k q^{k(k+1)/2}$ '' after the `` = '' sign.

p. 51, Ex. 2.13, first line: Replace `` $|\lambda/a|$ '' by
`` $|q \lambda/a|$ '' .

p. 167, right side of eq. (7.3.5): Replace `` $(q; bq;q)_m$ '' by
`` $(q, bq;q)_m$ ''  .

p. 210, (8.7.2): Replace `` ${(-z;q)_n}\over {(-N;q)_n}$ ''
by `` ${(q^{-z};q)_n}\over {(q^{-N};q)_n}$ '' .

p. 227, second line of Ex. 8.4 (ii): Replace `` $b q^{n+1-x-y}$ ''
by `` $b q^{N+1-x-y}$ '' .

p. 270, lines 11, 13, and 15: Replace `` R. '' by `` M. '' .

p. 116, (4.10.9): Replace `` $_B c_1$ '' by `` $b_B c_1$ '' in the
numerator of the fraction on the right side of the equals sign.

p. 82, middle: Replace `` Also, the case $d=q^{-n}$ '' by `` Also, 
the case $d=q^{-2n}$ '' .

p. 109, (4.5.1): Replace `` $;w)_\infty$ '' by `` $;q)_\infty$ '' .

p. 3, (1.2.11): Replace `` Re(c - a - b) '' by `` Re$(c - a - b)$ '' .

p. 29, Notes for Ex. 1.2:  Replace ``over a field with $q$ elements'' by
``over the field $GF(q)$ with $q$ a prime power''

p. 63, (3.3.1): In the first ${}_3\phi_2 $ replace `` $a,b,c,$ '' by 
`` $a,b,c$ '' .

p. 65, (3.4.3):  Replace `` ${}_{12}W_{11}[\ \ ] $ '' by `` 
${}_{12}W_{11}(\ \ ) $ '' .

p. 79, (3.8.5):  Replace each of the four `` $;q)_\infty$ '' by  `` 
$;q^2)_\infty$ '' .

p. 80, line 4 up:  Replace  `` ${bdf} \over {ac^2}$ '' by `` $f, {{bdf}\over {ac^2}}$ 
'' .

p. 81, line 5 up: Replace  `` $(a,b,cq/f;q)_\infty$ '' by `` 
$(a,b,cq/b;q)_\infty$ '' .

p. 82: Delete the period from the end of eq. (3.8.13) and add the phrase
``provided $d$ or $aq/d$ is not of the form $q^{-2n},\, n$ a 
nonnegative integer. ''

p. 82, (3.8.15): Replace the fraction in front of the  ${}_{10}W_9$ 
by the fraction
   $${{(acq, acq/df;q)_n \, (acq^{1-n}/d, acq^{1-n}/f;q^2)_n}\over
{(acq/d, acq/f;q)_n \, (acq^{1-n}, acq^{1-n}/df;q^2)_n}}\ .$$

p. 87, line 6 up: Replace ``Formula (3.10.2)'' by ``Formula (3.10.3)'' .  See
Gasper [1989c, (3.17)].

p. 123, line 4 of Ex. 4.11: Replace `` $bq/z;q)_\infty$ '' by `` 
$bq/a;q)_\infty$ '' .

p. 144, line 3 up: Replace `` [1985] '' by `` [1986b] '' .

p. 229, line 2: Delete `` $0 \le t < 1$ '' .

p. 229, line 12: Delete `` , Jain and Verma [1981] '' .

p. 182, (7.7.7): Replace `` $q^{-n/2}$ '' by ``  $q^{{1\over 2} -n}$  '' .

\vskip .1cm
Places which have too much (or too little) space (such as after the 
multiplication
dot in Exercise 2.9) are not included.
\vskip .1cm
This errata is available in .tex, .dvi, and .pdf file forms
via the World Wide Web at:

  http://www.math.northwestern.edu/$\scriptstyle\sim$george/preprints/
\vskip .2cm
We wish to thank everyone who took the time
to send errata. Please send any additional errata to:

\vskip .2cm

Professor George Gasper

Department of Mathematics

Northwestern University

Evanston, IL 60208-2730
\vskip .2cm
E-mail address: george@math.northwestern.edu

\vskip .4cm

{\bf Reviews of this book have appeared in:}
\vskip .2cm

Amer.\ Math.\ Monthly, 98 (1991), pp. 282--285 (by George E. Andrews).

Bull.\ London Math.\ Soc., 23 (1991), pp. 312--313  (by Tom H. Koornwinder).

Jahresbericht der Deutschen Mathematiker-Vereinigung (Buchbesprechungen),
  95 (1993), pp. 29--31 (in German by Clemens Markett).

Jour.\ Approx.\ Theory, 66 (1991), p. 352 (by Walter Van Assche).

Mathematical Reviews, 91d (1991), 91d:33034 (by Waleed A. Al-Salam).

Novye Knigi Za Rubejom (New Books Abroad), Vol. 7, Sept.\ 1990, pp. 30--32
(in Russian by N.M. Atakishiyev and S.K. Suslov).

SIAM Review, 33 (1991), pp. 489--493 (by Jet Wimp).

Zentralblatt f\"ur Math. 695 (1990), \#33001 (by R. K. Saxena).
\vskip .2cm
{\bf Russian translation:}
\vskip .2cm
MIR published a Russian translation (by N.M. Atakishiyev and S.K. Suslov)
of this book in 1993.
  \vskip .2cm
{\bf Errata for the Russian translation that are not covered by the
errata listed above:}
\vskip .2cm

RUSSp. 65: Replace the argument `` ${aq}\over{ef}$ '' of the first 
${}_8\phi_7$ in
eq. (2.10.3) by `` ${a^2q^{n+2}}\over{bcde}$ ''.

RUSSp. 69, line below eq. (2.10.19): Replace `` $cd=abdefgh$ '' by `` 
$cd=abefgh$ '' .

RUSSp. 293: In the `` ${}_8 \phi_7$ '' on the left side of
eq. (III.23),  replace `` $-aq/b$ '' by `` $aq/b$ '' .

\bye